\documentclass{amsart}

\author{P. W. Ng}

\address{Mathematisches Institut\\
Westfaelische Wilhelms-Universitaet Muenster\\
Einsteinstr. 62\\
48149 Muenster\\
Germany}

\title[Corona factorization property]{The corona factorization property}

\date{\today\\AMS Mathematics subject classification 46L35}

\newtheorem{thm}{Theorem}[section]
\newtheorem{df}{Definition}[section] 
\newtheorem{lem}[thm]{Lemma}
\newtheorem{cor}[thm]{Corollary}
\newtheorem{prop}[thm]{Proposition}
\theoremstyle{remark}

\newtheorem{qu}[thm]{Question} 

\newcommand{\J}{\mathcal{J}}
\newcommand{\B}{\mathcal{B}}

\newcommand{\A}{\mathcal{A}}

\newcommand{\C}{\mathcal C}
\newcommand{\E}{\mathcal E}

\newcommand{\D}{\mathcal D}
\newcommand{\M}{\mathbb M}

\newcommand{\K}{\mathcal K}
\newcommand{\Mul}{\mathcal M}

\def \bib(#1;#2;#3;#4;#5;#6)  {{#1},{\it #2} {#3},
{\bf#4} (#5) {#6}\par\smallskip}

\begin{document}

\begin{abstract}

     The corona factorization property is a property with connections to 
extension theory, $K$-theory and the structure of $C^*$-algebras.
This paper is a short survey of the subject, together with some 
new results and open questions.

\end{abstract}

\maketitle

\section{Introduction}  

     Absorbing extensions have played an important role 
in various places in operator algebra theory.  For example, they
were used by Kasparov to give a nice 
characterization of $KK$-theory (see
\cite{Black}).  They are also important for the stable existence 
and stable uniqueness results of classification theory (see
\cite{DE}).   Still most recently, they have been 
used in new ways to classify 
some new interesting classes of nonsimple 
$C^*$-algebras (one class being
purely infinite $C^*$-algebras with a unique ideal (see \cite{Restorff});
another class being Matsumoto algebras (these are
``generalized Cuntz-Krieger algebras") 
associated with primitive aperiodic substitutional subshifts (see
\cite{Ruiz}).  (The definition of absorbing extensions is in 
section 3.) 

     The corona factorization property was originally motivated by 
extension theory problems related to classification theory - specifically, 
the theory of absorbing extensions!  
A major starting point was Elliott and Kucerovsky's algebraic characterization
of nuclearly absorbing extensions (see \cite{EK}).  
Basicly, a separable nuclear stable $C^*$-algebra $\B$ 
has the corona factorization property if 
$\B$ has lots of absorbing extensions (see Theorem 3.1).  Among other things,
this includes automatically  
many of the extensions that have been important in 
classification theory, and leads to nice and useful characterizations
in $KK$-theory.      
Most recently, the corona factorization property has been used 
to classify
certain $C^*$-algebras associated with dynamical systems that were 
mentioned in the previous paragraph (see \cite{Ruiz}) 

  It turns out that the corona factorization property   
is also connected with basic questions about the structure of $C^*$-algebras.
Fundamental questions - like whether an extension of a stable $C^*$-algebra
by a stable $C^*$-algebra has a stable extension algebra
- seem to be related with the corona factorization property
and the techniques used to study it.  
This is not completely surprising, as 
one of the main techniques of this subject (related with Elliott and 
Kucerovsky's work on absorbing
extensions) is the theory of stability for  
$C^*$-algebras, as developped by (Larry) Brown, Cuntz, Hjelmborg,
Rordam and many others.  
Also, this subject often has the flavour of the theory of purely infinite
simple $C^*$-algebras.  

     This paper is a survey of the theory, with proofs of several new
results and many open questions.  We would like to emphasize that the  
corona factorization property is a type of ``regularity" property which
demarcates between ``nice" $C^*$-algebras (with ``nice" $K$-theory, 
``nice" structure theory and possibly ``nice" classification theory)
and not so ``nice" $C^*$-algebras.   
Indeed, our main examples of simple $C^*$-algebras with the 
corona factorization
property are examples which (with the additional assumption 
of nuclearity etc.) have been amenable to the $K$-theoretic classification
program.  On the other hand, our main examples which don't have 
corona factorization are exotic $C^*$-algebras (with very bad
perforation etc.)
constructed by Rordam (using, among other things, ideas of Villadsen). 

      In the next section, we introduce the corona factorization 
property and give basic examples and characterizations.
In section 3, we state the equivalence of corona factorization with
statements about extension theory and $KK$-theory.  
In section 4 (which is the largest section), we state connections with
the structure theory of $C^*$-algebras.  This section will also contain
proofs of some new results as well as some open questions.

     Before ending this introduction, we fix one notation.
All throughout this paper, the symbol ``$=_{df}$" will roughly mean
``is defined as" or ``is defined to be".  For example, 
the statement ``Let $\B =_{df} \A_0 \otimes \K$." reads as
``Let $\B$ be defined as the stabilization of $\A_0$.". 
 
     The survey component of this paper is joint work with
Dan Kucerovsky.

\section{Basic results and examples}

     The simplest statement of the corona factorization property (which,
for the purposes of this paper, we take as the definition) is the following
$K$-theoretic definition.  Firstly, a positive
element $c$ of a $C^*$-algebra
$\C$ is said to be \emph{norm-full} in $\C$ if the $C^*$-algebra ideal
that $c$ generates is all of $\C$  ; i.e., the $C^*$-subalgebra generated 
by $\C c \C$ is all of $\C$.  (In much $C^*$-algebra literature, the 
terminology 
``full" is often used instead of ``norm-full", but here we are 
dealing with multiplier algebras which have the strict topology as 
well as the norm topology and we need to distinguish norm-full
elements from strictly-full elements.) 

\begin{df}  Let $\B$ be a separable stable $C^*$-algebra.
Then $\B$ is said to have the \emph{corona factorization property}
if every norm-full projection in $\Mul(\B)$ is Murray-von Neumann 
equivalent to $1_{\Mul(\B)}$.
\end{df}
     
    With respect to this definition, one 
should think of the corona factorization property in the same
terms as (though not too closely with!) ``nice" properties like 
comparison of projections (which implies the corona factorization property!).
We will see that this property is flexible enough to include many 
interesting $C^*$-algebras but also rigid enough to rule out many
pathological examples.
   
    Many simple unital separable nuclear $C^*$-algebras, which have been
successfully classified using $K$-theory data, have (after stabilization)
the corona factorization property.  For example, we have 
the following proposition: 

\begin{prop} Let $\A_0$ be a unital separable simple $C^*$-algebra.
\begin{enumerate}
\item If $\A_0$ is exact, real rank zero, stable rank one, and has weakly
unperforated ordered 
$K_0$-group, then $\A_0 \otimes \K$ has the corona factorization
property.
\item If $\A_0$ is purely infinite, then $\A_0 \otimes \K$ has the 
corona factorization property.
\end{enumerate}
\end{prop}

\begin{proof}
 
     First, we prove (1).   
Let $\B =_{df} \A_0 \otimes \K$.  (For the notation ``$=_{df}$", see
the remark in the second last paragraph of the introduction.) 
Let $P$ be a norm-full projection in 
$\Mul(\B)$.   Since $\B$ has real rank zero, let $\{ p_n \}_{n=1}^{\infty}$
be a sequence of pairwise orthogonal projections in $\B$ (actually, in 
$P \B P$) such that 
$P = \sum_{n=1}^{\infty} p_n$, where the series converges in the strict 
topology on $\Mul(\B)$.   For each $N$, let 
$P_N =_{df} \sum_{n=1}^N p_n$.  We view each $P_N$ as a continuous 
function on $T(\A_0)$ (the simplex of unital traces on $\A_0$) in the 
natural way.   
Since $\B$ has stable rank one and weak unperforation and since 
$P$ is norm full in $\Mul(\B)$, we must have that the sequence
$\{ P_N \}_{N=1}^{\infty}$ increases to infinity pointwise on 
$T(\A_0)$.  Hence, by Dini's Theorem, $\{ P_N \}_{N=1}^{\infty}$ must
increase to infinity uniformly on $T(\A_0)$.  Hence, since $\B$ has 
real rank zero, stable rank one and weak unperforation,
$1_{\Mul(\B)}$ must be Murray-von Neumann equivalent
to a subprojection of $P$.  Hence, since $\B$ is stable,
$P$ is Murray-von Neumann equivalent to the unit of $\Mul(\B)$.

     Now we prove (2).  Again let $\B =_{df} \A_0 \otimes \K$ (where 
$\A_0$ is now purely infinite and simple).  Let $\pi : \Mul(\B) \rightarrow
\Mul(\B)/ \B$ be the natural quotient map.  Suppose that $P$ is a 
norm-full projection in $P$.  Then $\pi(P)$ is a nonzero projection in 
$\Mul(\B)/\B$.  But by \cite{Ror0}, 
$\Mul(\B)/\B$ is a simple, purely infinite
(nonseparable) $C^*$-algebra.
Hence, let $x \in \Mul(\B)$ be such that 
$\pi(x) \pi(P) \pi(x)^* = 1_{\Mul(\B)/ \B}$.
Hence, let $k \in \B$ be such that 
$x P x^* = 1_{\Mul(\B)} + k$.
Now since $\B$ is stable, let $S$ be an isometry in $\Mul(\B)$ such that 
$S^*k S$ is within $\epsilon$ of zero.
Then $S^* x P x^* S$ is within $\epsilon$ of $1_{\Mul(\B)}$.
Hence, since $\epsilon > 0$ is arbitrary, 
we can find $y \in \Mul(\B)$ such that $y P y^* = 1_{\Mul(\B)}$.
So the unit is Murray-von Neumann equivalent to a subprojection of $P$.
So since $\B$ is stable, $P$ is Murray-von Neumann equivalent to 
the unit of $\Mul(\B)$.
\end{proof}
 
     In addition, type I $C^*$-algebras with a ``finite-dimensionality"
condition also have (after stabilization) the corona factorization
property.
The following result can be found in \cite{KN1}: 

\begin{thm}  If $\B$ is a separable stable type I $C^*$-algebra with 
finite decomposition rank then 
$\B$ has the corona factorization property.
\end{thm}

     If $\B$ has continuous trace then the decomposition rank of $\B$ is 
the (ordinary topological) dimension of its spectrum.
Hence,  if $X$ is a finite-dimensional compact separable metric space
then $C(X) \otimes \K$ has the corona factorization property.
On the other hand, if $X$ is not finite-dimensional, then this 
can fail.  The following result can also be found in \cite{KN1}: 

\begin{thm}  Let $X$ be the countably infinite Cartesian product of 
spheres.  Then $C(X) \otimes \K$ does not have the corona factorization
property.
\end{thm}

     The proof of the previous result uses ideas of Rordam (and Villadsen!).
In particular, it is motivated by Rordam's example of an extension of 
a stable $C^*$-algebra by a stable $C^*$-algebra such that the extension
algebra is not stable (see \cite{RorExt}).  Using another of 
Rordam's ``exotic"
constructions, we also have a simple example.  The following can be 
found in \cite{KN3} (this is related to Rordam's discovery that 
``stability is not a stable property".  See also \cite{RorStab}):

\begin{thm}  There is a simple stable $AH$-algebra $\B$ with stable rank one
such that $\B$ does not have the corona factorization    
property.
\end{thm}

   We note that the corona factorization property does not rule out
every possible type of ``exotic" behaviour.  For example, in \cite{KNPerf},  
it is shown that there exists a simple $AH$-algebra with perforation
as well as the corona factorization property.

  Finally, we end this section with some basic alternate characterizations.
The full proof of the following can be found in \cite{KN3}: 

\begin{thm}  Let $\B$ be a separable stable $C^*$-algebra.
Then the following are equivalent:
\begin{enumerate}
\item $\B$ has the corona factorization property.
\item If $P$ is a norm-full projection in $\Mul(\B)/\B$ then there 
exists  $z \in \Mul(\B)/\B$ such that $z P z^* = 1_{\Mul(\B)/\B}$.
\item If $P$ is a norm-full projection in $\Mul(\B)/\B$ then
$P$ is properly infinite.
\item If $c$ is a positive element of $\Mul(\B)/ \B$ such that $C^*(c)$ does
not nontrivially intersect any proper ideal of $\Mul(\B)/\B$, 
then $c$ is properly infinite. 
\end{enumerate}
\end{thm} 
 
\begin{proof}[Sketch of proof]
 
     We prove (2) implies (1) and (3) implies (1).

     First for (2) implies (1):   Let $\pi : \Mul(\B) \rightarrow \Mul(\B)/\B$
be the natural quotient map.  Since $P$ is norm-full in $\Mul(\B)$,
$\pi(P)$ is norm-full in $\Mul(\B)/\B$.  Hence, by (2), let $x \in \Mul(\B)$
be such that $\pi(x) \pi(P) \pi(x)^* = 1_{\Mul(\B)/\B}$. 
Hence, let $k \in \B$ be such that $x P x^* = 1_{\Mul(\B)} + k$.
Since $\B$ is stable, let $S$ be an isometry in $\Mul(\B)$ such that 
$S^* k S$ is within $\epsilon$ of zero.
Hence, $S^* x P x^* S$ is within $\epsilon$ of the unit of $\Mul(\B)$.
Since $\epsilon > 0$ is arbitrary and since $\B$ is stable,
$P$ is Murray-von Neumann equivalent to the unit of $\Mul(\B)$.

    Next for (3) implies (1):  Suppose that $P$ is a norm-full projection
in $\Mul(\B)$.  Since $P$ is norm-full in $\Mul(\B)$, there is a positive
integer $N$ such that the $N$-times direct sum of $P$ with itself is 
Murray-von Neumann equivalent to the unit of $\Mul(\B)$.  But by (3),
$P$ is properly infinite.  Hence, $P$ is Murray-von Neumann equivalent to the
unit of $\Mul(\B)$. 
\end{proof}
 
   Finally, we note that property (2) in Theorem 2.5 is formally 
similar to one of the characterizations of purely infiniteness for
simple unital separable $C^*$-algebras.  Hence, it is not surprising
that the theory of corona factorization is related to the structure
and stability theory of $C^*$-algebras and often has the flavour of purely
infinite simple $C^*$-algebras. 

\section{The corona factorization property and extension theory}

The corona factorization property first arose in our study of 
$C^*$-algebras which are ``nice" from the point of view of 
absorbing extensions.  This is, for instance, one of the 
first places where the corona factorization property enters into
the classification program for simple unital separable nuclear
$C^*$-algebras. For basic extension theory (Busby invariant, 
BDF-sums etc...) we refer the reader to \cite{Black} and \cite{WO}. 

We say that an extension $\tau$ is \emph{absorbing} if 
$\tau +_{BDF} \rho$ is BDF-equivalent to $\rho$ for every trivial extension
$\rho$, where $+_{BDF}$ is the
BDF-sum.   Note that the definition implies that $\tau$ is not unital.
For unital $\tau$, the definition is exactly the same except that 
the trivial extension $\rho$ is required to lift to a unital $*$-homomorphism
(into the multiplier algebra $\Mul(\B)$).   
Finally, if in the above definition we require that $\rho$ be weakly 
nuclear then we say that $\tau$ is \emph{nuclearly absorbing}.  (Note that
if either the ideal algebra or quotient algebra is a nuclear $C^*$-algebra,
then every extension is weakly nuclear.)
Absorbing extensions play 
an important role in a number of places - among other things, Kasparov used 
them to  
give a clean characterization of $KK$-theory (see \cite{Black}).  
They are also important in the 
stable existence and stable 
uniqueness theorems of classification theory (see, for example, 
\cite{DE}), as well as recent classification results for 
nonsimple $C^*$-algebras (see \cite{Restorff} and \cite{Ruiz}).

   To see the connection between the corona factorization property and
absorbing extensions, we need to first discuss norm-full extensions.

\begin{df}  Suppose that $\A$ and $\B$ are separable $C^*$-algebras such
that $\B$ is stable.  An extension $\tau : \A \rightarrow \Mul(\B) / \B$ is 
said to be \emph{norm-full} if for every nonzero, positive element 
$a \in \A$, $\tau(a)$ is a norm-full element of $\Mul(\B)/\B$.
\end{df}

   We note that the Lin and Kasparov extensions as well as many other
useful extensions (e.g. useful to classification theory; see \cite{DE}) 
are norm full extensions. 

   Note that an absorbing extension is necessarily norm-full (in the definition
of absorbing, we can always choose the trivial extension to be norm-full).
On the other hand, it is not always the case that the converse is true (see
below).  The converse, however, is always true (for nuclearly absorbing)
exactly in the case where the ideal algebra has the corona factorization
property. In \cite{KN4},  we prove the following:

\begin{thm}  Suppose that $\B$ is a separable, stable $C^*$-algebra.
Then the following conditions are equivalent:
\begin{enumerate}
\item $\B$ has the corona factorization property.
\item Every norm-full extension of $\B$ is nuclearly absorbing.
\item Every norm-full extension of $\B$ is nuclearly absorbing in the sense
of approximate unitary equivalence (with unitaries in the corona).
\item Every norm-full trivial extension of $\B$ is nuclearly absorbing.
\end{enumerate}
\end{thm}

    Note that from Theorems 2.3 and 3.1, it follows that if $X$ is the 
countably infinite Cartesian product of spheres then $C(X) \otimes \K$ has
a nonabsorbing full extension.  For more details, see \cite{KN1}. 
 
   In connection to the above, the corona factorization property also
gives a clean characterization of $KK^1$ as well as a nice uniqueness
theorem.   

     When the context is clear, we use the terminologies ``full extension"
and ``full $*$-homomorphism" instead of 
``norm-full extension"  and  ``norm-full $*$-homomorphism". 
(When confusion with the strict topology on the multiplier algebra
is not possible...)

\begin{thm}  Let $\B$ be a separable stable nuclear $C^*$-algebra.
Then the following conditions are equivalent:

\begin{enumerate}
\item $\B$ has the corona factorization property.
\item For every separable nuclear $C^*$-algebra $\A$,
$KK^1(\A, \B)$ is the group of full extensions under unitary equivalence
by multiplier unitaries.
\item Let $\A$ be a unital separable nuclear $C^*$-algebra.
Let $\phi, \psi : \A \rightarrow \Mul(\B)/ \B$ be two unital
full $*$-monomorphisms  
Then $[ \phi ] = [\psi]$ in $KL( \A, \Mul(\B)/ \B)$ if and only if
$\phi$ and $\psi$ are approximately unitarily equivalent, with unitaries
coming from the corona algebra.
\end{enumerate}

\end{thm}

     We finally note that Theorems 3.1 and  3.2 have played 
a role in recent classification theory for certain nonsimple 
$C^*$-algebras with a distinguished ideal
(see, for example,  \cite{Restorff} and \cite{Ruiz}). 

\section{The corona factorization property and the structure 
of $C^*$-algebras.}

    It turns out that the corona factorization property is also of interest
from the point of view of the structure of $C^*$-algebras.  This (long)
section is an exposition of results so far obtained.  There will also be
some
proofs of new results as well as some open questions.

       The corona factorization property demarcates among ``nice" and 
not so ``nice" $C^*$-algebras.      
Firstly, Rordam has constructed an example of a simple $AH$-algebra
$\B_0$,
with stable rank one, such that there is a positive integer $n \geq 2$ with
$\M_n (\B_0)$ stable but $\B_0$ itself is not stable (see \cite{RorStab}).   
In the context of the corona factorization property, this phenomenon is
ruled out.  In \cite{KN3}, we prove the following:

\begin{thm}  Suppose that $\B$ is a separable, stable $C^*$-algebra.
Then the following are equivalent:
\begin{enumerate}
\item $\B$ has the corona factorization property.
\item Suppose that $\D$ is a full, hereditary subalgebra of $\B$.
Suppose that there is an integer $n \geq 1$ such that $\M_n (\D)$ is 
stable.  Then $\D$ itself is stable.
\end{enumerate}
\label{thm:stab}
\end{thm}

     Hence, in the context of the corona factorization property, stability
is a stable property for full hereditary subalgebras.

     Next, Rordam has constructed an example of an extension of a separable
stable $C^*$-algebra, by a separable stable $C^*$-algebra, such that 
the extension algebra is not stable (see \cite{RorExt}).
When the ideal algebra has the corona factorization property, this 
behaviour is also ruled out.  In \cite{KN3}, we prove the following:

\begin{thm}  Suppose that $\J, \E$ and $\A$ are separable 
$C^*$-algebras, such that $\J \otimes \K$ has the corona 
factorization property.
Suppose that we have an exact sequence of the form
\[
0 \rightarrow \J \rightarrow \E \rightarrow \A \rightarrow 0.
\]
Then $\E$ is stable if-and-only-if $\J$ and $\A$ are stable.
\label{thm:ext}
\end{thm}
    
     Moreover, under appropriate hypotheses, we actually have a converse. 
In \cite{KN3}, we also prove the following:  

\begin{thm}  Suppose that $\B$ is a stable, separable, simple, 
real rank zero $C^*$-algebra with cancellation of projections.
Then the following are equivalent:
\begin{enumerate}
\item $\B$ has the corona factorization property.
\item Every extension of $\B$, by a separable stable $C^*$-algebra, gives
a stable extension algebra.
\end{enumerate}
\end{thm}   

    By the same argument as the above, we also have the following:

\begin{thm}  Suppose that $\B$ is a stable, separable, simple, 
real rank zero $C^*$-algebra with cancellation of projections.
Then the following are equivalent:
\begin{enumerate}
\item $\B$ has the corona factorization property.
\item Every extension of $\B$, by the compact operators $\K$, gives a
stable extension algebra.
\item Suppose that 
\[
0 \rightarrow \B \rightarrow \C \rightarrow \K \rightarrow 0
\]
is an essential quasidiagonal extension such that $\C$ has real rank zero.
Then $\C$ is stable.
\end{enumerate}
\end{thm} 
 
    Motivated by the above results, we look for connections between
the corona factorization property and interesting properties about the 
structure of simple separable unital $C^*$-algebras.  One such property
is Rordam's notion of \emph{regularity}:

\begin{df}  A $C^*$-algebra $\B$ is said to be \emph{regular} if every 
full hereditary subalgebra of 
$\B$, with no nonzero unital quotients and no 
nonzero bounded traces, is stable.
\end{df}
 
    Firstly, regularity implies the corona factorization property:

\begin{lem}  Suppose that $\B$ is a separable, stable $C^*$-algebra.
If $\B$ is regular, then $\B$ has the corona factorization property.
\label{lem:RegCFP}
\end{lem}
\begin{proof}
Suppose that $\D$ is a full, hereditary subalgebra of $\B$ such that 
there is a positive integer $n \geq 1$ with 
$\M_n (\D)$ being a stable $C^*$-algebra.
Then $\D$ has no nonzero unital quotient and no nonzero bounded trace
(for otherwise, $\M_n(\D)$ would have such - which is impossible since
$\M_n (\D)$ is stable).  
Hence, by regularity $\D$ must be stable.
But $\D$ was arbitrary.
Hence, by Theorem 4.1, $\B$ has the corona factorization property.
\end{proof}

   The converse of the above result is not known (even with the additional
assumptions of simplicity and real rank zero).  

  Regularity is an interesting property from the point of view of $C^*$-algebra
structure.   
The next lemma is due to Hjelmborg and Rordam.  We present
the short argument for the convenience of the reader.

\begin{lem}\  Suppose that $C$ is a simple, stable, exact,
separable $C^{*}$-algebra which is regular. 
Then $C$ is either purely infinite or stably finite.
\label{lem:RegDich}
\end{lem}
\begin{proof}
Suppose that $C$ is not stably finite.  Then no nonzero hereditary
subalgebra
of $C$ is stably finite.
Let $D$ be a nonunital hereditary subalgebra of $C$.
Since $D$ is not stably finite, and since $D$ is exact,
$D$ has no nonzero bounded traces.  Hence, by hypothesis,
$D$ is stable.
But $D$ was arbitrary.
Hence, every nonzero hereditary subalgebra of $C$ is either unital or
stable.
Hence, by a result of Zhang's \cite{zhang} (see also
\cite{kirchberg}), the algebra $C$
is either the compact operators or purely infinite.
\end{proof}

     Henceforth, we will say that a $C^*$-algebra $C$ has
\emph{dichotomy} if $C$ 
is either stably finite or purely infinite (in the sense of 
Kirchberg and Rordam).  Note that this property passes on
to full, hereditary subalgebras.  

    Next, we weaken the notion of regularity by introducing the notion
of asymptotic regularity.

\begin{df}  Suppose that $\B$ is a separable, stable $C^*$-algebra.
Then $\B$ is said to be \emph{asymptotically regular}, if whenever 
$\D$ is a full hereditary subalgebra of $\B$ with no nonzero unital quotients
and no nonzero bounded traces, there is a positive integer $n \geq 1$ such 
that $\M_n (\D)$ is stable.
\end{df}

     It is interesting to try and understand the relationship between 
regularity and asymptotic regularity.   In the case of exact,
simple,
stable, separable, real rank
zero algebras, 
the two notions are the same.  To prove this, we require
the following proposition, which can be found in \cite{Brown} 
Theorem 4.23:

\begin{prop}  Let $B$ be a separable, stable $C^*$-algebra.  Let $P$
be a projection in $\Mul(B)$.  
Then $P B P$ is a stable, full, hereditary subalgebra of $B$ if-and-only-if
$P$ is Murray-von Neumann equivalent to the unit of $\Mul(B)$.
\label{prop:unit}
\end{prop}

We also need the following lemma which is in \cite{Kuc} Lemma 10:

\begin{lem}
Let $\B$ be a separable stable $C^*$-algebra.
Then for every hereditary subalgebra $\B_0$ of $\B$, there exists
a multiplier projection $P \in \Mul(\B)$ such that 
$P \B P \cong \B_0$. 
\label{lem:Kuc}
\end{lem}

\begin{prop}  Let $B$ be a separable, exact, simple,
stable, real rank zero $C^*$-algebra.
Then $B$ is regular if-and-only-if $B$ is asymptotically regular.
\label{prop:RegAsym}
\end{prop}
\begin{proof}
The ``only-if" direction is clear.

We proceed to prove the ``if" direction.
If $B$ is type I, then $B$ is the compact operators over a separable,
infinite dimensional Hilbert space.  Hence, $B$ is automatically regular.
Hence, we may assume that $B$ is not type I.
First, suppose that there is a projection $p$ of $B$ such that 
the hereditary subalgebra $p B p$ has a unital trace (hence, $B$ is 
stably finite).

Suppose that $B$ is asymptotically regular.
Suppose, to the contrary, that $B$ is not regular.
Hence, let $D$ be a nonzero hereditary subalgebra of $B$, with no unit,
and no nonzero bounded traces, such that $D$ is not stable.
By Lemma \ref{lem:Kuc}, let $P$ be a projection in $\Mul(B)$ such that 
$P B P$ is isomorphic to $D$.  
By asymptotic regularity of $B$, let $n$ be the least integer such
that $\M_n (D)$ is stable ($n \geq 2$).
For each positive integer $m$,
let $Q_m =_{df} \bigoplus_{i=1}^m P$.  Since $B$ is stable, $\Mul(B)$ contains
a copy of $O_m$ for each $m$; and hence 
for each positive integer $m$, $Q_m$ is defined, up to Murray-von Neumann
equivalence, as a projection in $\Mul(B)$.  Hence, by Proposition 
\ref{prop:unit},
$n$ is the least integer such that $Q_n$ is Murray-von Neumann equivalent
to $1_{\Mul(B)}$, the unit of $\Mul(B)$.  In particular, $Q_1 = P$ is not
Murray-von Neumann equivalent to $1_{\Mul(B)}$.  

     Now since $B$ has real rank zero, let $\{ p_k \}_{k=1}^{\infty}$.
be a sequence of pairwise orthogonal projections in $B$ such that 
$P = \sum_{k=1}^{\infty} p_k$, where the sum on the right hand side
converges in the strict topology in $\Mul(B)$.
Note that
for each unital
trace $\tau$ in $p B p$, $\tau$ extends naturally to 
a semifinite trace on $B$.  This trace, in turn, gives a trace on
the positive cone of $\Mul(B)$.  
We use ``$\tau$" to denote any one of these traces.
Now since $Q_n$ is Murray-von Neumann equivalent to $1_{\Mul(B)}$,
we must have that for each unital trace $\tau$ of $pBp$,
$\tau(P) = \infty$.  Hence, for each unital trace $\tau$ of $pBp$,
$\sum_{k=1}^{\infty} \tau( p_k ) = \infty$.

     Now by \cite{PR} Theorem 5.8, and since $B$ is not type I,
$B$ is weakly divisible.  Hence, for each $k$, for each strictly positive
integer $l \geq 2$, let $m_{k,l}$ be an integer greater than $l$, and  
let $q_{k,l,1}, q_{k,l,2}, ..., q_{k,l, m_{k,l}}$ be pairwise orthogonal 
projections in $B$ such that 
\begin{enumerate}
\item for each $i$, $q_{k,l,i}$ is a subprojection of $p_k$;
\item for each $i,j$, $q_{k,l,i}$ is Murray-von Neumann equivalent to 
$q_{k,l,j}$;
\item $q_{k,l, 1} + q_{k,l, 2} + ... + q_{k,l, m_{k,l} -1}$ 
is a subprojection of $p_k$; and
\item $p_k$ is a subprojection of 
$q_{k,l,1} + q_{k,l,2} + ... + q_{k, l, m_{k,l}}$.
\end{enumerate}

     Now let $T( p B p)$ be the simplex of unital traces of $p B p$.  
$T( p B p)$ is a compact, convex set (compact in the pointwise topology).
As we have observed in previous paragraphs, for each $\tau \in T(p B p)$,
$\tau (P) = \sum_{k=1}^{\infty} \tau( p_k ) = \infty$.
Moreover, by the compactness of $T(p B p)$, $\sum_{k = 1}^{\infty} p_k$
is a sequence of positive, (affine) continuous functions on 
$T( p B p)$ which blows up to infinity uniformly on $T(p B p)$.  
We recursively define a sequence $\{ r_i \}_{i=2}^{\infty}$ of projections
in $B$, and an increasing
sequence $\{ n_i \}_{i=2}^{\infty}$ of positive integers as follows:
    
     Let $n_2$ be an integer such that 
$\sum_{k=1}^{n_2} \tau( q_{k, 1, 1}) \geq 1$, for all $\tau \in T(p B p)$. 
Let $r_1$ be the projection $r_1 =_{df} \sum_{k=1}^{n_2} q_{k, 1, 1}$.

      Suppose that both $n_i$ and $r_i$ have been chosen.
Let $n_{i+1}$ be an integer such that 
$\sum_{k= n_i + 1}^{n_{i + 1}} \tau( q_{k, i + 1, 1}) \geq 1$, for 
all $\tau \in T( p B p)$. 
Let $r_{i + 1}$ be the projection
$r_{i+1} =_{df} \sum_{k = n_i + 1}^{ n_{i+1}} q_{k, i +1 , 1}$.

      The recursion is complete. 
Now let $R =_{df} \sum_{i = 1}^{\infty} r_i$.  The sum converges in the
strict topology in $\Mul( B )$.  Hence, $R$ is a projection in $\Mul(B)$.
By our choice of the $r_i$s, $\tau(R) = \infty$ for all 
$\tau \in T(p B p)$ (where, as usual, we use ``$\tau$" to also denote
the induced trace).  Hence, $R B R$ is a full, hereditary subalgebra of 
$B$, with no unit and no nonzero bounded traces.

     We want to show that there is no positive
integer $m$ such that $\M_m (RBR)$ 
is stable.  
So suppose, to the contrary, that such an $m$ does indeed exist.
Then, by Proposition \ref{prop:unit}, 
$\bigoplus_{j=1}^m R$ is Murray-von Neumann
equivalent to the unit of $\Mul(B)$.  
But by our definition of $R$,
$\bigoplus_{j=1}^m R$
is a subprojection of a projection of the form 
$P \oplus r$, where $r$ is a projection in $B$.
Hence, by \cite{WO} Lemma 16.2, $P \oplus r$ is Murray-von Neumann
equivalent to the unit of $\Mul(B)$.  Hence, let $x \in \Mul(B)$ be 
such that $x( P \oplus r ) x^* = 1_{\Mul(B)}$.  Hence,
$xPx^* + xrx^* = 1_{\Mul( B)}$.  Since $x r x^* \in B$, and since
$B$ is stable, let $S$ be an isometry in $\Mul(B)$ such that 
$S^* x r x^* S$ has norm strictly less than $\epsilon$.  
Hence, $S^* x P x^* S$ is within $\epsilon$ of $1_{\Mul(B)}$.  
Since we can choose $\epsilon$ to be arbitrarily small, there is an
element $y \in \Mul(B)$ such that $y P y^* = 1_{\Mul(B)}$.  Hence,
$1_{\Mul(B)}$ is Murray-von Neumann equivalent to a subprojection of
$1_{\Mul(B)}$.  Hence, by \cite{WO} Lemma 16.2, 
$P$ is Murray-von Neumann equivalent
to $1_{\Mul(B)}$.  This is a contradiction.
Hence, $B$ must be regular.

     For the case where there is no projection $p \in B$ such that 
$p B p$ has a unital trace ($B$ is not stably finite), we can 
simply choose $r_i$ to be $r_i =_{df} q_{i, i, 1}$ for every $i$.
The rest of the proof is exactly the same.
\end{proof}

    For the nonsimple, nonreal rank zero case, 
the situation is not yet known.  However, we do have an equivalence
if we add the corona factorization property to asymptotic regularity.

\begin{prop}
Let $B_0$ be a unital, separable $C^*$-algebra, and let
$B =_{df} B_0 \otimes \K$ be the stabilization of $B_0$.
Then $B$ is regular if-and-only-if $B$ is asymptotically regular and has 
the corona factorization property.
\label{prop:Regular}
\end{prop}
\begin{proof}

     The ``only-if" direction follows by definition. 

     Now for the ``if" direction.  Suppose that $D$ is a full, hereditary
subalgebra of $B$, with no nonzero unital quotients and no nonzero bounded
traces.  By asymptotic regularity, let $n$ be a positive integer such
that $\M_n (D)$ is stable.  Then since $B$ has the corona factorization
property, and by Theorem 4.1, $D$ must be stable.  So since $D$ is arbitrary,
$B$ is regular. 
\end{proof}

      Now let $B$ be a stable $C^{*}$-algebra, and let $\pi$ be a nonzero
$*$-representation of $B$.  Then $\pi$ extends to a unique strictly
continuous, surjective $*$-homomorphism
$\pi^{\prime\prime}:\Mul(B)\rightarrow\Mul(\pi (B))$ (between multiplier
algebras).
Let $\J_{\pi}$ be the proper, norm-closed ideal of $\Mul(B)$ given
by $\J_{\pi}=_{df}\{c\in\Mul(B):\pi^{\prime\prime} (c)\in  \pi(B)\}$.

    Now let $B_0$ be a separable, unital $C^{*}$-algebra, and let
$B:=B_0\otimes\K$ be the stabilization of $B_0$.  Suppose
that
$\tau$ is a unital trace on $B_0$.  Then $\tau$ extends canonically
to a
trace on the positive cone of 
$\Mul(B)$, which we also denote by ``$\tau$".
We let $\J_{\tau}$ be the proper ideal of $\Mul(B)$ given by the
norm-closure of $\{b\in\Mul(B):\tau (b^{*}b)<\infty \}$.

\begin{df}\ Let $B_0$ be a separable, unital $C^{*}$-algebra, and let
$B=_{df}B_0\otimes\K$ be the stabilization of $B_0$.
Let $\J$ be a proper ideal of $\Mul(B)$, which contains $B$.
we say that the ideal $\J$ of $\Mul(B)$ is \emph{regular} if
\begin{enumerate}
\item$\J$ is contained in an ideal of the form $\J_{\pi}$, for
some nonzero $*$-representation $\pi$ of $B$; or
\item$\J$ is contained in an ideal of the form $\J_{\tau}$, for
some unital trace $\tau$ on $B_0$.
\end{enumerate}
Otherwise, we say that $\J$ is \emph{nonregular}.  
\end{df}

     If $B_0$ is the Cuntz algebra, then both $B_0 \otimes \K$ and
$\K$ have simple corona algebras (see \cite{Lin1}). Hence, it follows, by
definition, that both $\Mul(B_0\otimes\K)$ and $\Mul(\K)$ have no
nonregular ideals.  If $B_0$ is a simple, unital $AF$-algebra, such
that
the tracial simplex of $B_0$ has only $n<\infty$ extreme points, then
$\Mul(B_0\otimes\K)$ has $n$ maximal, proper ideals.  These ideals
are all regular, since they come from traces on $B_0$ (see
\cite{Lin1}).

\begin{df}\  Let $B_0$ be a unital $C^{*}$-algebra, and let
$B=_{df}B_0\otimes\K$ be the stabilization of $B_0$.
We say that $B$ has \emph{property $R$}, if whenever $p$ is a
projection
contained inside a proper ideal of $\Mul(B)$, $p$ is also contained
inside a regular ideal.
\end{df}

     Under the conditions of exactness, simplicity and
real rank zero, property $R$ is
rather strong.   

\begin{prop}  Let $B_0$ be a simple, exact,
unital, separable, real rank
zero $C^*$-algebra, and let $B =_{df} B_0 \otimes \K$ 
be the stabilization of $B_0$.
Suppose that $B$ has property $R$.
Then
\begin{enumerate}
\item $B$ has dichotomy; and  
\item $B$ has the corona factorization property.
\end{enumerate}
\label{prop:R}
\end{prop}

\begin{proof}

     We first prove statement i).  
Suppose, to the contrary, that $B$ is neither purely infinite, nor
stably finite.  By \cite{Ror0} Theorem 3.2, $\Mul( B )$ has a proper
ideal $\J$, which properly contains $B$.  
By \cite{Zh2}, since $B$ has real rank zero, $\J$ is the norm linear
span of its projections.  Hence, let $P$ be a projection in $\J$ such
that $P$ is not contained in $B$.  Since $B_0$ is exact and is not stably
finite, $B_0$ has no unital traces.  From this, and the simplicity of $B$,
$\Mul(B)$ has no regular ideals other than $B$.  Hence,
$P$ is a projection, contained in a proper ideal of $\Mul(B)$, such that
$P$ is not contained in a regular ideal.  This contradicts property $R$.
Hence, $B$ must be either purely infinite or stably finite.

     We now proceed to prove the second statement, ii).  
Suppose, to the contrary, that $B$ does not have the corona factorization
property.  By i), we know that $B$ must be stably finite (since, by 
Proposition 2.1, 
a stable, separable, simple, purely infinite $C^*$-algebra always has 
the corona factorization property).   Hence, by exactness,
$B_0$ has a unital trace.  So let $T(B_0)$ be the simplex of unital traces
on $B_0$. $T(B_0)$ is a compact, convex set (compact with respect to the
pointwise topology).

     Since $B$ does not have the corona factorization property, let $P$
be a norm-full projection in $\Mul(B)$ such that $P$ is not Murray-von 
Neumann equivalent to $1_{\Mul(B)}$, the unit of $\Mul(B)$.
By norm-fullness, there is an integer $n \geq 2$ such that 
$\bigoplus_{i=1}^n P$ is Murray-von Neumann equivalent to
$1_{\Mul(B)}$.  
Hence, for every $\tau \in T(B_0)$, $\tau(P) = \infty$  (again, we 
are using ``$\tau$" to denote also the natural
extension to the positive cone of $\Mul(B)$).

     Since $B$ has real rank zero, let $\{ p_k \}_{k=1}^{\infty}$ be 
a sequence of pairwise orthogonal projections in $B$ such that 
$P = \sum_{k=1}^{\infty} p_k$, where the sum converges in the strict
topology in $\Mul(B)$.  Hence, for each $\tau \in T(B_0)$, 
$\sum_{k=1}^{\infty} \tau( p_k ) = \infty$.

     If $B$ is type I, then it is the compact operators on a separable
Hilbert space;  and hence, $B$ automatically has the corona factorization
property.  Hence, let us assume that $B$ is not type I.  Then, by
\cite{PR} Theorem 5.8, $B$ is weakly divisible.
Hence, for each $k$, for each strictly positive integer $l \geq 2$, let
$m_{k,l}$ be a positive integer greater than $l$, and let
$q_{k, l, 1}, q_{k, l, 2}, ..., q_{k, l, m_{k,l}}$ be pairwise orthogonal
projections in $B$ such that 

\begin{enumerate}
\item for each $i$, $q_{k, l, i}$ is a subprojection of $p_k$;
\item for each $i, j$, $q_{k, l, i}$ is Murray-von Neumann equivalent
to $q_{k, l, j}$; 
\item $q_{k, l, 1} + q_{k, l, 2} + ... + q_{k, l, m_{k,l} - 1}$ is a
subprojection of $p_k$; and
\item $p_k$ is a subprojection of 
$q_{k, l, 1} + q_{k, l, 2} + ... + q_{ k, l, m_{k,l}}$
\end{enumerate}

     Now for each $i$, define $r_i$ as in the proof of Proposition
\ref{prop:RegAsym}. 
And, as in the proof of Proposition  \ref{prop:RegAsym}, 
let $R =_{df} \sum_{i=1}^{\infty} r_i$, where
the sum converges in the strict topology in $\Mul(B)$.
As in Proposition \ref{prop:RegAsym}, 
for each integer $m$, $\bigoplus_{j = 1}^m R$ is
not Murray-von Neumann equivalent to $1_{\Mul (B)}$.  Hence,
$R$ is contained in a proper ideal of $\Mul(B)$.  But also, as in
Propositon \ref{prop:RegAsym}, 
$\tau(R) = \infty$ for every $\tau \in T(B_0)$.  
Hence, $R$ is a projection which is contained in a proper ideal of 
$\Mul(B)$, but which is not contained in a regular ideal.
This contradicts property $R$.
\end{proof}

     Next, we relate property $R$ to regularity, for the case 
of simple, exact,
separable, stable $C^*$-algebras.

\begin{prop}\  Let $B_0$ be a simple, exact,
unital, separable $C^{*}$-algebra,
and
let $B:=B_0\otimes\K$ be the stabilization of $B_0$.
Then $B$ is regular if-and-only-if $B$ has property $R$ and the corona
factorization property.
\end{prop}

\begin{proof}

    We first prove the ``only-if" direction. 
Suppose that $B$ is regular.   
By Lemma \ref{lem:RegCFP}, 
$B$ has the corona factorization property. 
Hence, it suffices to prove that $B$ has property $R$.
So let $P$ be a projection in $\Mul(B)$ such that $P$ is 
contained in a proper ideal of $\Mul(B)$.
If $P$ is an element of $B$, then automatically, $P$ is contained in
a regular ideal of $\Mul(B)$.  Hence, we may assume that $P$ is
not contained in $B$.

     Now suppose, to the contrary, that $P$ is not contained inside any
regular ideal of $\Mul(B)$.  
Hence, $P B P$ must be a full hereditary subalgebra of $B$, with no
unit and no nonzero bounded traces.  Hence, since $B$ is regular,
$P B P$ is stable.  Hence, by Proposition \ref{prop:unit}, 
$P$ is Murray-von 
Neumann equivalent to the unit of $\Mul( B )$.  This contradicts
our assumption that $P$ is contained in a proper ideal of $\Mul(B)$.
Hence, $P$ is contained in a regular ideal of $\Mul(B)$.  By the
arbitrariness of $P$, $B$ has property $R$.

     Next, we prove the ``if" direction.  Suppose that $B$ has both
property $R$ and the corona factorization property.
Suppose that $D$ is a nonzero hereditary subalgbra of $B$, with no unit and
no nonzero bounded traces.  By Lemma \ref{lem:Kuc},  
let $P$ be a projection in
$\Mul(B)$ such that $P B P$ is isomorphic to $D$.  Since $D$ has no
nonzero bounded traces, $P$ cannot be contained in any regular 
ideal of $\Mul(B)$.  Hence, since $B$ has property $R$, 
$P$ is a norm-full element of $\Mul(B)$.  Hence, by the corona
factorization property and by \cite{WO} Lemma 16.2, $P$ must be
Murray-von Neumann equivalent to the unit of 
$\Mul(B)$.  Hence, $P B P$ is stable.
Hence, $D$ is stable.  Since $D$ is arbitrary, $B$ must be regular.
\end{proof}

\begin{cor}  Suppose that $B$ is a simple, separable, stable, 
exact, real rank
zero $C^*$-algebra.
Then $B$ is regular if-and-only-if $B$ has property $R$.
\end{cor}

   Collecting the above results and restricting to the simple, real rank
zero case, we have the following:

\begin{thm}
Suppose that $\B$ is a simple, separable, stable, exact, real rank zero
$C^*$-algebra. 
Then the following are equivalent:
\begin{enumerate}
\item $\B$ is regular. 
\item $\B$ is asymptotically regular.
\item $\B$ has property $R$.
\end{enumerate}
\end{thm}
 
\begin{qu}  Suppose that $\B$ is a separable, stable $C^*$-algebra.
If $\B$ has the corona factorization property, then is $\B$ regular?
What if, in addition, $\B$ is simple and has real rank zero?
\end{qu}

\begin{qu}
With possibly additional assumptions, what is (are)
the connection(s) of the
above statements with the statement that every extension of $\B$, by a 
separable stable $C^*$-algebra, gives a stable extension algebra? 
What about dichotomy?
\end{qu}

\begin{qu}
Does every simple separable stable nuclear real rank zero $C^*$-algebra
have the corona factorization property?  What if we also assumed that
the $C^*$-algebra was $AH$?
\end{qu}

\end{document}